\documentclass[11pt]{article}
\usepackage{bidipoem}
\usepackage{graphicx}
\usepackage{amsfonts, amsmath, amsthm, amssymb, latexsym, amscd}
\usepackage{pstricks}

\newtheorem{thm}{Theorem}[section]
\newtheorem{exam}[thm]{Example}
\newtheorem{lem}[thm]{Lemma}
\newtheorem{corol}[thm]{Corollary}

\newenvironment{prf}{{\bf \noindent Proof } \rm }{\hfill{$\Box$}}

\title{A row generation method for the general case of inverse continuous facility location problem}

\author{
Jafar Fathali
       \thanks{Faculty of Mathematical Sciences, Shahrood University of Technology, University Blvd., Shahrood, Iran, E-mail: fathali@shahroodut.ac.ir.}\index{Fathali, Jafar}
}

\begin{document}
\date{}
\maketitle

\begin{abstract}
In a single facility location problem, a set of points is given and the goal is finding the optimal location of new facility respect to given criteria such as minimizing time, cost and distances between the clients and facilities. On the other side, the inverse models try to modify the parameters of the problem with the minimum cost such that a given point becomes optimal. In this paper, we introduce a novel algorithm for the general case of the inverse single facility location problem with variable weights in the plane. The convergence and optimality conditions of the algorithm are presented. Then in the special cases the inverse minisum and minimax single facility location problems are considered and the algorithm tested on some instances. The results indicate the efficiency of the algorithm on these instances.    
\end{abstract}


\section{Introduction}
The single facility location problem is well known in operations research literature and plays an important role in practice. The continuous version of this problem asks for finding the location of a facility in the plane, such that the cost of servicing clients by the facility is minimized. In the classical version of facility location problems, a set of points, represent the location of clients, are given and a weight is assigned to each point. Two basic objective functions in the facility location problems are minimizing the sum of weighted distances and the maximum weighted distances between clients and the facility, which called minisum and minimax single facility location problems, respectively. These problems also called median and center problems in the literature of location theory. For more details in continuous facility location models, the reader is referred to \cite{LMV88}.

In some real applications, the facilities may already exist and the goal is  modifying the parameters of the problem with minimum cost, such that the given facility location become optimal. These kind of location models are called inverse facility location problems. 

There are many researches in the graph version of the inverse facility location problems. Among them, Cai et al. \cite{CYZ99} showed that the inverse center problem is NP-hard. Burkard et al. \cite{BPZ04} proposed an $O(nlogn)$ algorithm for the inverse 1-median problem on a tree. Then Galavii \cite{G10} improved the time complexity of this problem to linear time. Alizadeh and Burkard \cite{AB11} showed that the inverse 1-center problem can be solved in $O(n^2)$ time. Guan and Zhang \cite{GZ12} and Nguyen and Sepasian \cite{NS15} considered the inverse 1-median and 1-center problems on trees with Chebyshev and Hamming norms, respectively. Sepasian and Rahbarnia \cite{SR15} solved the inverse 1-median problem with varying vertex and edge length on trees. When the underlying network is a block graph, Nguyan \cite{N16-2} proposed a solution algorithm for the inverse 1-median problem with variable vertex weight. Nazari et al. \cite{NFNV18} investigated the inverse of backup 2-median problem with variable edge length and vertex weight. Recently, Omidi et al. \cite{OF20} proposed an algorithm for solving the inverse balanced facility location models on a tree. 

The continuous version of inverse facility location problem is rarely considered by authors. The minisum single facility location problem has been considered by Burkard et al. \cite{BGG10}. They showed the Euclidean case of this problem with variable vertex weights can be solved in $O(n log n)$ time. Baroughi-Bonab et al. \cite{BBA10} investigated the inverse minisum single facility location problem with variable coordinates. They showed the problem with rectilinear and Chebyshev norms are NP-hard. Nazari and Fathali \cite{NF18} applied a meta-heuristic algorithm for solving the inverse backup 2-median problem with variable coordinates.

In this paper, we investigate the general case of inverse continuous single facility location problem with variable weights. This problem asks to modify the weights of existing points with minimum cost such that a given point becomes optimal. A novel row generation algorithm is introduced. Then, we show the algorithm converges to the solution of inverse facility location problem.  The algorithm explained for the inverse minimax and minsum location models.   

In what follows, the general form of continuous single facility location problem and its inverse are presented in Section \ref{sec2}. In Section \ref{sec3}, a general algorithm is introduced and its convergence is investigated. Then solving two special cases of continuous inverse location models by the presented algorithm are explained in Subsections \ref{sec3.1} and \ref{sec3.2}. Computational results are presented in Section \ref{sec4}. Section \ref{sec5} contains the summary and the conclusion on this paper.

\section{The general form of continuous facility location problem and its inverse}\label{sec2}
Let $\mathbf{A}$ be the set of $n$ points $\mathbf{A}_1, \mathbf{A}_2, ..., \mathbf{A}_n$ in the plane. Where $\mathbf{A}_i=(a_i,b_i)$ for $i=1,...,n$ are the location of clients. For $i=1,...,n$, let the point $\mathbf{A}_i$ has the weight $w_i$. The goal in a facility location problem is finding a point $\mathbf{x}=(x,y)$ in the plane, which is the location of a new facility such that the function $f(\mathbf{x},\mathbf{w},\mathbf{A})$ is minimized, i.e.
\begin{equation*}\label{minimax}
(P) \ \min_{\mathbf{x}\in \mathbb{R}^2} f(\mathbf{x},\mathbf{w},\mathbf{A}),
\end{equation*}
where $ \mathbf{w}=(w_1,...,w_n)$ is the vector of weights of existing points. 

In the inverse model the aim is modifying some parameters of the problem with minimum cost such that a given point $\bar{\mathbf{x}}$ become optimal. The changing parameters can be the weights, coordinates of points or both. 

In this paper, we consider the weight modifying case. Thus, we want to modify the weights of existing points $w_i$, for $i=1,...,n$, to ${\hat{w_i}}=w_i+p_{i}-q_{i}$ with minimum cost.
%, such that the circle $C^*$ be the optimal with comparing to any other circle with radius $ r_0 $.
Where $p_{i}\geq 0$ and $q_{i}\geq 0$ are the values of augmenting and reduction of the weight of point $\mathbf{A}_i$ and  bounded from above by $u^+_i$ and $u_i^-$, respectively. Suppose that $c_i^+$ and $c_i^-$ denote the costs of augmenting and reduction of per unit of $w_i$, respectively. Then the inverse problem with variable weights can be modeled as follow.
\begin{align}
(P_{w})\ \min & \sum_{i=1}^n(c_i^+p_{i}+c_i^-q_{i})\\
&f(\bar{\mathbf{x}},\hat{\mathbf{w}},\mathbf{A})\leq f(\mathbf{x},\hat{\mathbf{w}},\mathbf{A}) \ \ \ \ \ \forall\mathbf{x}\in\mathbb{R}^2\label{inf1}\\
&\hat{w_i}=w_i+p_{i}-q_{i} \ \ \ \ i=1,\cdots, n,\label{what}\\
& 0\leq p_{i} \leq u_i^+\qquad\qquad i=1,\cdots , n,\\
& 0\leq q_{i} \leq u_i^-\qquad\qquad i=1,\cdots , n,\\
&\hat{w_{i}}\geq 0\qquad\qquad i=1,...,n.
\end{align}
This model is a problem with infinite constraints and difficult to solve.
% Note that the constraints \ref{what} can be omitted by replacing constraint \ref{inf1} with the following,
%$$f(\bar{\mathbf{x}},\mathbf{w}+{\mathbf{p}}-\mathbf{q},\mathbf{A})\leq f(\mathbf{x},\mathbf{w}+{\mathbf{p}}-\mathbf{q},\mathbf{A})  \ \ \ \ \ \forall\mathbf{x}\in\mathbb{R}^2,$$
%where $\mathbf{p}=(p_1,...,p_n)$ and $\mathbf{q}=(q_1,...,q_n).$

In the next section we present a novel row generation algorithm for solving the inverse single facility location problem.
 
\section{The general algorithm for the inverse problem}\label{sec3}
As mentioned in Section \ref{sec2}, the general form of inverse model contains infinite constraints. In this section we present a method with finite constraints for solving the inverse problem. 

Let $\mathbf{x^{(0)}}$ be the optimal solution of problem $(P)$. The weights of existing points should be modified such that the given point $\bar{\mathbf{x}}$ become better than any other point. The point $\mathbf{x^{(0)}}$ is optimal with respect to the current wights. Therefore, in the first step we modify the weights such that $\bar{\mathbf{x}}$ become better than $\mathbf{x^{(0)}}$. Thus we consider a new model by replacing the infinite constraints (\ref{inf1}) with the following single constraint,
$$ f(\bar{\mathbf{x}},\mathbf{\hat{w}},\mathbf{A})\leq f(\mathbf{x^{(0)}},\mathbf{\hat{w}},\mathbf{A}).$$
Then the new weights of points are obtained by solving this new model. In the next step, $\mathbf{x^{(1)}}$ the optimal solution of problem $(P)$ with respect to the new weights, is calculated. Then the following new constraint is added to the model  
$$ f(\bar{\mathbf{x}},\mathbf{\hat{w}},\mathbf{A})\leq f(\mathbf{x^{(1)}},\mathbf{\hat{w}},\mathbf{A}), $$
and the new weights are obtained by solving this model. This procedure continues until the weights do not changed. % or the algorithm reaches the maximum number of iterations.

For $i=1,...,n$, let $p_{i}^{(k)}$ and $q_{i}^{(k)}$ be the values of modifying the weights of existing points in iteration $k$ then 
$$C^{(k)}=\sum_{i=1}^{n}(c_i^+ p_{i}^{(k)}+c_i^- q_{i}^{(k)})$$
is the value of objective function for the current solution.

These ideas lead us to the following algorithm.
 
\bigskip
\noindent {\bf Algorithm [A1]}\label{al1}.\\
%\textbf{Input:} The star graph $T$, where $x$ as the location of the facility, is the vertex  of $T$ with degree larger than 1. Also, the cost of increasing and decreasing of edges lengths are given.\\
%\textbf{Output:} The new lengths of edges $\hat{l_i}$ such that $f_2=0$.\\
%\textbf{Initialization}:\\ 
\begin{enumerate}
\item \textbf{Find} $\mathbf{x^{(0)}}$, the optimal solution of problem $(P)$.
\item \textbf{Set} $k=1$ and $\mathbf{w^{(0)}}=\mathbf{w}$.
\item \textbf{If} $f(\mathbf{\bar{x}},{\mathbf{w}^{(0)}},\mathbf{A})=f(\mathbf{x}^{(0)},{\mathbf{w}^{(0)}},\mathbf{A})$ stop, the current solution is optimal. 

%\textbf{Iteration step} \\
%\item \textbf{Set} $\mathbf{\hat{w}}=\mathbf{w}$.
\item \textbf{Find} $(\hat{w}_i,p_{i}^{(1)} ,q_{i}^{(1)})$, $i=1,...,n$, the optimal solution of the following problem,
\begin{align}
(P_0)\ \min & \sum_{i=1}^n(c_i^+p_{i}+c_i^-q_{i})\\
&f(\bar{\mathbf{x}},\hat{\mathbf{w}},\mathbf{A})\leq f(\mathbf{x^{(0)}},\hat{\mathbf{w}},\mathbf{A}) \\
%&f(\bar{\mathbf{x}},\mathbf{w}^{(0)}+{\mathbf{p}}-\mathbf{q},\mathbf{A})\leq f(\mathbf{x^{(0)}},\mathbf{w}^{(0)}+{\mathbf{p}}-\mathbf{q},\mathbf{A}) \\
&\hat{w_i}-p_{i}+q_{i}=w_i^{(0)} \ \ \ \ i=1,\cdots, n,\label{wk}\\
& 0\leq p_{i} \leq u_i^+\qquad\qquad i=1,\cdots , n,\\
& 0\leq q_{i} \leq u_i^-\qquad\qquad i=1,\cdots , n,\\
&\hat{w_{i}}\geq 0\qquad\qquad\qquad i=1,...,n.
%& \mathbf{\hat{w}}\geq 0.
\end{align} 
If this problem is infeasible stop, the main problem is infeasible.
\item \textbf{Set} $\mathbf{w}^{(1)}=\mathbf{w}^{(0)}+{\mathbf{p}^{(1)}}-\mathbf{q}^{(1)}$.
%\item \textbf{For} $i=1,..,n$ set $u_i^+=u_i^+-p_i^{(1)}+q_i^{(1)}$ and $u_i^-=u_i^- -q_i^{(1)}+p_i^{(1)}$.
%\begin{enumerate}
\item \textbf{While} $||\mathbf{w}^{(k)}- \mathbf{w}^{(k-1)}||> \epsilon$ do
\begin{enumerate}
\item \textbf{Find} $\mathbf{x^{(k)}}$ the optimal solution of the following sub-problem
$$\min_{\mathbf{x}\in \mathbb{R}^2} f({\mathbf{x}},{\mathbf{w}^{(k)}},\mathbf{A}).$$ 
\item %\textbf{Let} $P_k$ be the problem $P_{k-1}$ with replacing $w_i^{k-1}$ by $w_i^{k}$ in the constraints \ref{wk} and adding the following constraint
\textbf{Add} the following constraint to problem $(P_{k-1})$
$$f(\bar{\mathbf{x}},\hat{\mathbf{w}},\mathbf{A})\leq f(\mathbf{x}^{(k)},\hat{\mathbf{w}},\mathbf{A}),$$ 
%$$f(\bar{\mathbf{x}},\mathbf{w}^{(k)}+{\mathbf{p}}-\mathbf{q},\mathbf{A})\leq f(\mathbf{x}^{(k)},\mathbf{w}^{(k)}+{\mathbf{p}}-\mathbf{q},\mathbf{A}),$$ 
%and replace $w_i^{(k-1)}$ by $w_i^{(k)}$ in the constraints \ref{wk} 
and call the new problem $(P_{k})$.
\item \textbf{If} problem $(P_k)$ is infeasible stop, the main problem is infeasible.
\item \textbf{Find} $(w_i^{(k+1)},p_{i}^{(k)} ,q_{i}^{(k)})$, $i=1,...,n$, the optimal solution of the problem $(P_k)$.
%\item \textbf{Set} $\mathbf{w}^{(k+1)}=\hat{\mathbf{w}}.$ %\mathbf{w}^{(k)}+{\mathbf{p}^{(k)}}-\mathbf{q}^{(k)}$.
%\item \textbf{For} $i=1,..,n$ set $u_i^+=u_i^+-p_i^{(k)}+q_i^{(k)}$ and $u_i^-=u_i^- -q_i^{(k)}+p_i^{(k)}$. %$w_i^{(k+1)}=\hat{w_i} \ \ \ \ i=1,\cdots, n$
\item \textbf{Set} $k=k+1$.
%\end{enumerate}
\end{enumerate}
\textbf{End while}
\item \textbf{For} $i=1,...,n$ set
$$(w_i^*,p_i^*,q_i^*)=(w_i^{(k+1)},p_{i}^{(k)} ,q_{i}^{(k)}). $$
\item\textbf{Set} 
$$C^*=\sum_{i=1}^{n}(c_i^+ p_{i}^*+c_i^- q_{i}^*). $$
\end{enumerate}
\textbf{End of algorithm}
\bigskip

Note that in each iteration, one constraint has been added to the previous model and the best weights are chosen, thus the algorithm obey the row generation rule. Moreover, instead of solving a problem with infinite constraints, the algorithm solve some problems with finite constraints. Indeed, the problem is reduced to the discrete version of the inverse location problem. 

The following theorem shows that if in two iterations the weights of points do not changes then an optimal solution is found.     
%The following algorithm 
\begin{thm}\label{stopw}
Let $t$ be an iteration of Algorithm A1 such that $\mathbf{w}^{(t+1)}= \mathbf{w}^{(t)}$, then $\mathbf{\hat{w}}=\mathbf{w^{(t+1)}}$ is the optimal solution of $(P_w)$.
\end{thm}
\begin{prf}
In the iteration $t$ of the Algorithm A1, the following problem is solved
\begin{align*}
(P_t)\ \min & \sum_{i=1}^n(c_i^+p_{i}+c_i^-q_{i})\\
&f(\bar{\mathbf{x}},\hat{\mathbf{w}},\mathbf{A})\leq f(\mathbf{x^{(k)}},\hat{\mathbf{w}},\mathbf{A})\qquad k=0,...,t \\
%&f(\bar{\mathbf{x}},\mathbf{w}^{(0)}+{\mathbf{p}}-\mathbf{q},\mathbf{A})\leq f(\mathbf{x^{(0)}},\mathbf{w}^{(0)}+{\mathbf{p}}-\mathbf{q},\mathbf{A}) \\
&\hat{w_i}=w_i^{(0)}+p_{i}-q_{i} \ \ \ \ i=1,\cdots, n,\\
& 0\leq p_{i} \leq u_i^+\qquad\qquad i=1,\cdots , n,\\
& 0\leq q_{i} \leq u_i^-\qquad\qquad i=1,\cdots , n,\\
&\hat{w_{i}}\geq 0\qquad\qquad\qquad i=1,...,n.
%& \mathbf{\hat{w}}\geq 0.
\end{align*} 
Where $\mathbf{x^{(k)}}$ is the solution of $k$-th sub-problem. Specially, for $k=t$,
$$f({\mathbf{x}^{(t)}},{\mathbf{w}^{(t)}},\mathbf{A})=\min_{\mathbf{x}\in \mathbb{R}^2}  f(\mathbf{x},{\mathbf{w}^{(t)}},\mathbf{A}).$$ 
By the algorithm $(\mathbf{w^{(t+1)}},\mathbf{p^{(t)}},\mathbf{q^{(t)}})$ is the optimal solution of $(P_t)$, then 
$$f(\bar{\mathbf{x}},{\mathbf{w}^{(t+1)}},\mathbf{A})\leq f(\mathbf{x}^{(t)},{\mathbf{w}^{(t+1)}},\mathbf{A}),$$
$$\mathbf{w^{(t+1)}}=\mathbf{w^{(0)}}+\mathbf{p^{(t)}}-\mathbf{q^{(t)}}. $$
Since $\mathbf{w}^{(t)}=\mathbf{w}^{(t+1)}$ then  
$$f(\mathbf{x}^{(t)},{\mathbf{w}^{(t+1)}},\mathbf{A})=\min_{\mathbf{x}\in \mathbb{R}^2}  f(\mathbf{x},{\mathbf{w}^{(t+1)}},\mathbf{A}).$$ 
Therefore,
$$f(\bar{\mathbf{x}},{\mathbf{w}^{(t+1)}},\mathbf{A})\leq f(\mathbf{x},{\mathbf{w}^{(t+1)}},\mathbf{A})\ \ \forall\mathbf{x}\in\mathbb{R}^2,$$ 
It means $(\mathbf{w^{(t+1)}},\mathbf{p^{(t)}},\mathbf{q^{(t)}})$ is a feasible solution of $(P_w)$.

Now let $(\mathbf{w'},\mathbf{p'},\mathbf{q'})$ be a feasible solution of $(P_w)$ then 
$$f(\bar{\mathbf{x}},{\mathbf{w}'},\mathbf{A})\leq f(\mathbf{x},{\mathbf{w}'},\mathbf{A}) \ \ \forall\mathbf{x}\in\mathbb{R}^2,$$
$$\mathbf{w'}=\mathbf{w^{(0)}}+\mathbf{p'}-\mathbf{q'}. $$ 
This implies that $(\mathbf{w'},\mathbf{p'},\mathbf{q'})$ is also feasible for $(P_t)$. Since  $(\mathbf{w^{(t+1)}},\mathbf{p^{(t)}},\mathbf{q^{(t)}})$ is the optimal solution of $(P_t)$ then 
$$\sum_{i=1}^n(c_i^+p^{(t)}_{i}+c_i^-q^{(t)}_{i})\leq \sum_{i=1}^n(c_i^+p'_{i}+c_i^-q'_{i}).$$
It means $(\mathbf{w^{(t+1)}},\mathbf{p^{(t)}},\mathbf{q^{(t)}})$ is the optimal solution of $(P_w)$.
\end{prf}

The stopping condition of the algorithm is based on Theorem \ref{stopw}. This stopping condition can be replaced by another one which is the case that the value of objective function in sub-problem is equal to $f(\mathbf{\bar{x}},{\mathbf{w}^{(t)}},\mathbf{A})$.

\begin{thm}\label{stopfx}
If in the iteration $t$ of Algorithm A1, $f(\mathbf{\bar{x}},{\mathbf{w}^{(t)}},\mathbf{A})=f(\mathbf{x}^{(t)},{\mathbf{w}^{(t)}},\mathbf{A})$, then $\mathbf{w}^{(t)}$ is an optimal solution of $(P_w)$.
\end{thm}
\begin{prf}
Note that $\mathbf{x}^{(t)}$ is the optimal solution of sub-problem, thus
$$f(\mathbf{x}^{(t)},{\mathbf{w}^{(t)}},\mathbf{A})=\min_{\mathbf{x}\in \mathbb{R}^2} f(\mathbf{x},{\mathbf{w}^{(t)}},\mathbf{A}).$$ 
Since  $f(\mathbf{\bar{x}},{\mathbf{w}^{(t)}},\mathbf{A})=f(\mathbf{x}^{(t)},{\mathbf{w}^{(t)}},\mathbf{A})$, therefore 
$$f(\bar{\mathbf{x}},{\mathbf{w}^{(t)}},\mathbf{A})\leq f(\mathbf{x},{\mathbf{w}^{(t)}},\mathbf{A})\ \ \forall\mathbf{x}\in\mathbb{R}^2,$$ 
so $\mathbf{w}^{(t)}$ is a feasible solution for problem $(P_w)$.
The proof of optimality is the same as proof of Theorem \ref{stopw}. 
\end{prf} 

Moreover, if $\mathbf{\bar{x}}= \mathbf{x}^{(t)}$ then $f(\mathbf{\bar{x}},{\mathbf{w}^{(t)}},\mathbf{A})=f(\mathbf{x}^{(t)},{\mathbf{w}^{(t)}},\mathbf{A})$. Thus the following corollary is immediately concluded.  

\begin{corol}\label{stopx}
If in the iteration $t$ of Algorithm A1, $\mathbf{\bar{x}}= \mathbf{x}^{(t)}$, then $\mathbf{w}^{(t)}$ is optimal for $(P_w)$ .
\end{corol}

 The following theorem indicates the  convergence of the algorithm. 

%\begin{thm}
%Let $\mathbf{w}$ be a feasible solution of problem $P_w$ which obtained by algorithm A1. If the problem $P$ has a unic optimal solution then $\mathbf{w}$ is the optimal solution of $P_w$. 
%
%\end{thm}

\begin{thm}
Algorithm A1 converges to an optimal solution of $(P_w)$.  
\end{thm}
\begin{prf}
Since in each iteration of Algorithm A1, a constraint is added to make $\mathbf{\bar{x}}$ better than some new points in the plane. Therefore, the feasibility improved. Thus if $k\rightarrow \infty$ then the algorithm will terminate with finding the optimal solution. 
\end{prf}

%Note that since problem $(P_w)$ is a problem with infinite constraints, therefore finding a feasible solution for this problem is not easy. 
In some special cases such as inverse minisum problem with Euclidean norm, we know that the optimal solution of problem $(P)$ lies in the convex hull of the existing points (see e.g. \cite{LMV88}). Therefore in this case, if the given point $\bar{x}$ does not lie in the convex hull, then the inverse model is infeasible. The following lemma shows the relation between infeasibility of problem $(P_w)$ and problem $(P_k)$ in Algorithm A1.   

\begin{lem}
In any iteration $k$ of Algorithm A1, if problem $(P_k)$ is infeasible then problem $(P_w)$ is infeasible, too.  
\end{lem} 

Note that the time complexity of the algorithm depends on sub-problem $(P)$ and problem $(P_k)$, in iteration $k$. One of inverse location models that can be solved in polynomial time by Algorithm A1, is inverse minisum location problem.  In the next sections we explain the presented algorithm for solving the inverse minisum and minimax location problems.

\subsection{The inverse minisum location problem}\label{sec3.1}
In this section we consider the special case that the objective function of location problem is minimum the weighted sum of distances between existing points and the new facility, i.e. 
\begin{equation*} %\label{minisum}
(P_{ms}) \ \min_{\mathbf{x}\in \mathbb{R}^2} \sum_{i=1}^n w_id(\mathbf{x},{A_i}),
\end{equation*}
where $d(x,A)$ is the distance between the points $x$ and $A$. This problem called Fermat-weber problem.

In this case, the inverse model with variable weights can be written as follow.
\begin{align}
(P_{ws})\ \min & \sum_{i=1}^n(c_i^+p_{i}+c_i^-q_{i})\\
&\sum_{i=1}^{n}\hat{w_i} d(\bar{x},A_i)\leq \sum_{i=1}^{n}\hat{w_i} d(x,A_i) \ \ \ \ \ \forall\mathbf{x}\in\mathbb{R}^2\label{infms}\\
&\hat{w_i}=w_i+p_{i}-q_{i} \ \ \ \ i=1,\cdots, n,\\
& 0\leq p_{i} \leq u_i^+\qquad\qquad i=1,\cdots , n,\\
& 0\leq q_{i} \leq u_i^-\qquad\qquad i=1,\cdots , n.
\end{align} 
This is a linear programming model with infinite constraints. Although, in the case that distances are measured by Euclidean norm, Burkard et al. \cite{BGG10} presented a linear model with finite constraints, however, in the general $L_p$ norm there isn't any efficient method. Therefore, our proposed algorithm is the first for solving the inverse minisum location problem with $L_p$ norm.   

Using Algorithm A1, in the  iteration $k$, the following sub-problem should be solved.
\begin{equation}\label{minisum-sub}
\min_{\mathbf{x}\in \mathbb{R}^2} \sum_{i=1}^n w_i^{(k)}d(\mathbf{x},{A_i}),
\end{equation}
When the distances are measured by $L_p$ norms, this sub-problem can be solved by the iterative  Weiszfeld algorithm (see e.g. \cite{LMV88}).

Also, the following constraint should be added to problem $(P_{k-1})$ in  iteration $k$,
%\begin{equation}\label{minisum-const}
%\sum_{i=1}^n \hat{w_i}d(\mathbf{\bar{x}},{A_i})\leq \sum_{i=1}^n \hat{w_i}d(\mathbf{x}^{(k)},{A_i}),
%\end{equation} 
%where $\hat{w_i}=w_i^{(k)}+p_{i}-q_{i}$ for $i=1,\cdots, n$. By substitution $\hat{w_i}$ to \ref{minisum-const}, the following linear constraint is obtained,
%or
\begin{equation}\label{minisum-const1}
\sum_{i=1}^n \hat{w_{i}}(d(\mathbf{\bar{x}},{A_i})-d(\mathbf{x}^{(k)},{A_i}))\leq 0,
\end{equation} 
which is a linear constraint. Since the linear programs and sub-problems can be solved in a polynomial time, therefore the time complexity of Algorithm A1 for the inverse minisum location problem with variable weights, is polynomial. 

\begin{exam}\label{ex1}
Consider the instance of the inverse minisum location problem given in Table \ref{inst1}, which is the only instance that solved in \cite{BGG10}. The given point is $\bar{\mathbf{x}}=(0,0)$. Burkard et al. \cite{BGG10} found that the optimal solution is $\mathbf{w}^*=(0,5,5,\frac{10}{\sqrt{2}})$. We examined our algorithm for finding the solution of this instance. 

\begin{table}\caption{ The parameters of considered instance in Example \ref{ex1}}\label{inst1}
\centering
%\small %\scalebox{1}{
\begin{tabular}{ccccccc}
\hline
i& $A_i=(a_i,b_i)$& $w_i$& $u^-_i$&$u^+_i$&$c^-_i$&$c^+_i$\\
\hline
1& $(1,0)$&0&0&5&$\sqrt{2}$&$\sqrt{2}$\\
2& $(\frac{1}{\sqrt{2}},\frac{1}{\sqrt{2}})$&0&0&5&7&7\\ %  \\
3& $ (-\frac{1}{\sqrt{2}},\frac{1}{\sqrt{2}})$&0&0&5&1&1\\
4& $(0,-1)$&$ \frac{10}{\sqrt{2}}$&0&0&   0&0\\
\hline
\end{tabular}
\end{table}

The iteration results of Algorithm A1 are shown in Table \ref{result1}. In this table the column with heading $\mathbf{x}^{(k)}$ indicates the solution of sub-problem (\ref{minisum-sub}) in iteration $k$ with respect to $\mathbf{w}^{(k)}$. The algorithm find the optimal solution with tolerance $\epsilon=0.01$, after 11 iterations.  The results show $\mathbf{x}^{(k)}$ converges to the given point $\mathbf{\bar{x}}=(0,0)$. The cost of modifying $\mathbf{w}^{(0)}$ to $\mathbf{w}^{(11)}$ is $C^*=39.9836$.

\begin{table}\caption{ The iteration results of Algorithm A1 on instance of Example \ref{ex1}.}\label{result1}
\centering
\small %\scalebox{1}{
\begin{tabular}{ccccccc}
\hline
iteration& $\mathbf{x}^{(k)}=(x^{(k)},y^{(k)})$& $w^{(k)}_1$& $w^{(k)}_2$ &$w^{(k)}_3$&$w^{(k)}_4$&$||\mathbf{w}^{(k)}-\mathbf{w}^{(k-1)}||$\\
\hline
0&$(0.0000,-1.0000)$&5.0000&0.8979&0.0000&7.0711&7.1278\\
1& $(0.1697,-0.3725)$&0.8790&2.9114&5.0000&7.0711&4.5866\\
2& $(0.0206,-0.3787)$&2.1373&3.4609&5.0000&7.0711&1.3730\\ 
3& $ (0.0877,-0.1749)$&0.5318&3.9241&5.0000&7.0711&1.6710\\
4& $(-0.0125,-0.1498)$&1.0139&4.3648&5.0000&7.0711&0.6532\\
5& $(-0.0469,0.0698)$&0.3567&4.4951&5.0000&7.0711&0.6700\\
6& $(0.0008,-0.0541)$&0.3567&4.8101&5.0000&7.0711&0.3150\\ 
7& $ (0.0168,-0.0132)$&0.1549&4.8136&5.0000&7.0711&0.2019\\
8& $(0.0014,-0.0077)$&0.0000&4.9647&5.0000&7.0711&0.2164\\
9& $(-0.0015,-0.0015)$&0.0264&4.9724&5.0000&7.0711&0.0275\\
10& $(0.0004,-0.0011)$&0.0000&4.9918&5.0000&7.0711&0.0328\\ 
11& $(-0.0003,-0.0003)$&0.0041&4.9950&5.0000&7.0711&0.0052\\
\hline
\end{tabular}
\end{table}

\end{exam}

\subsection{The inverse minimax location problem}\label{sec3.2}
In this section, the main steps of  Algorithm A1 for the case that the objective function of location problem is minimum the weighted maximum distance from existing points to the new facility, is explained. The minimax location problem is as follow. 
\begin{equation*} %\label{minisum}
(P_{mm}) \ \min_{\mathbf{x}\in \mathbb{R}^2} \max_{i=1,...,n} w_id(\mathbf{x},{A_i}).
\end{equation*}
%where $d(x,A)$ is the distance between the points $x$ and $A$. 

Then the inverse model with variable weights can be written as follow.
\begin{align}
(P_{wm})& \min \sum_{i=1}^n(c_i^+p_{i}+c_i^-q_{i})\\
&\max_{i=1,...,n}\{\hat{w_i} d(\bar{x},A_i)\}\leq \max_{i=1,...,n}\{\hat{w_i} d(x,A_i)\} \ \ \ \ \ \forall\mathbf{x}\in\mathbb{R}^2\label{infmm}\\
&\hat{w_i}=w_i+p_{i}-q_{i} \ \ \ \ i=1,\cdots, n,\\
& 0\leq p_{i} \leq u_i^+\qquad\qquad i=1,\cdots , n,\\
& 0\leq q_{i} \leq u_i^-\qquad\qquad i=1,\cdots , n.
\end{align} 
This problem is a nonlinear model with infinite constraints. As far as we know, this is an open problem and there isn't any efficient method for solving it. 
%Cai et al. \cite{CYZ99} showed that this problem is NP-hard.
Using our proposed algorithm, some models with finite constraints should be solved. The following sub-problem should be solved in the iteration $k$ of Algorithm A1.
\begin{equation}\label{minimax-sub}
\min_{\mathbf{x}\in \mathbb{R}^2} \max_{i=1,...n} w_i^{(k)}d(\mathbf{x},{A_i}).
\end{equation}
This minimax sub-problem can be solved in a linear time (see, e.g. \cite{EM11} and references therein).

Also, in the iteration $k$, the following constraint should be added to  problem $(P_{k-1})$, 
\begin{equation}\label{minimax-const1}
\max_{i=1,...,n} \{\hat{w_{i}}(d(\mathbf{\bar{x}},{A_i})-d(\mathbf{x}^{(k)},{A_i}))\}\leq 0.
\end{equation}
Although the added constraints in the inverse minisum model are linear, these constraints in the inverse minimax case are nonlinear. However, the problem $(P_k)$ is a discrete version of the inverse minimax location problem which called inverse 1-center problem. Nguyen et al. \cite{NNNLP19} showed that the discrete inverse 1-center problem can be solved in $O(n^2 logn)$ time. Therefore, the time complexity of our proposed algorithm for inverse minimax location problem is polynomial and it is more efficient than solving a problem with infinite constrains.

\section{Computational results}\label{sec4}
In this section we examine some test problems on our presented algorithm. The algorithm  was  written  in  MATLAB 2014.  All  the  experiments were run on a PC with Intel Core i7 processor, 8 GB of RAM and CPU 2.4 GHz.
The proposed algorithm was  tested on 4  test problems with
varying given points and $L_p$ norms. In the Euclidean case the results are compared with those obtained by the linear programming method of Burkard et al. \cite{BGG10}. In the other cases, there isn't any efficient method in the literature, therefore just the results of our algorithm are presented.
In tables \ref{resd18} and \ref{resultl2}, the notation $C_B$ indicates the optimal costs which obtained by the method of Burkard et al. \cite{BGG10}. After finding the new weights by the method of  Burkard et al. \cite{BGG10}, to verify its feasibility, we applied the Wiszfeld method with the new weights. The obtained solution is shown by $\mathbf{{x}}_B$.
 
The results of all test problems are reported for $\epsilon=0.01$. The number of the last iteration is shown by $t$. 

As the first test problem, consider the points and their required data that are given in Table \ref{ta18}. The coordinate of these points are given from \cite{F09}.  

\begin{table}[h]\centering
\small
\begin{tabular}{cc|cc|c}\hline
$(a_i,b_i,w_i,c_i^-,c_i^+,u_i^-,u_i^+)$  & &$(a_i,b_i,w_i,c_i^-,c_i^+,u_i^-,u_i^+)$& &$(a_i,b_i,w_i,c_i^-,c_i^+,u_i^-,u_i^+)$\\\hline
$(1,2,3,2,1,3,5)$& &$(4,4,1,2,1,1,3)$& &$(7,1,2,1,1,2,6)$\\
$(1,3,2,2,3,2,7)$& &$(4,9,2,1,1,2,2)$& &$(7,2,3,3,5,3,7)$\\
$(2,5,1,3,4,1,9)$& &$(5,3,2,2,5,2,6)$& &$(8,5,5,2,3,4,8)$\\
$(3,6,3,2,1,3,8)$& &$(5,5,1,4,3,1,5)$& &$(8,8,3,4,2,3,3)$\\
$(4,8,2,2,2,2,4)$& &$(6,6,3,5,2,3,2)$& &$(9,7,4,2,5,4,9)$\\
$(4,1,3,1,3,2,8)$& &$(6,3,3,1,4,2,4)$& &$(9,6,5,5,5,5,1)$\\\hline
\end{tabular}
\caption{Data for 18 points problem.}\label{ta18}
\end{table}

The results for varying given points with Euclidean norm are shown in Table \ref{resd18}. Each column of this table, contains a given point $\mathbf{\bar{x}}$ and the obtained results for this point. The results indicates that in all cases, the algorithm find the optimal solution of problem $(P_w)$ with tolerance $\epsilon$. 

\begin{table}\caption{ The results of Algorithm A1 on the instance with 18 points.}\label{resd18}
\centering
\small %\scalebox{1}{
\begin{tabular}{|c||c|c|c|}
\hline
$\bar{\mathbf{x}}$&$(2,2)$&$(3,5)$&$(7,7)$\\
%\hline
$\mathbf{x}^{(t)}$&$(2.0978,2.0105)$&$(2.9972,4.9992)$&$(6.9831,6.9826)$ \\
$\mathbf{x}_{B}$&$(2.0681,1.9962)$&$(2.9954,4.9983)$&$(6.9950,6.9987)$\\
\hline
$C^*$&100.1798&72.2736&57.6731\\
%\hline
$C_B$&101.2458&72.7461&58.48071\\
\hline
Iteration(t)&17&14&17\\
\hline
$w^{(t)}_1$&8.0000&8.0000&0.0000\\
$w^{(t)}_2$&1.3723&2.0000&0.0000\\
$w^{(t)}_3$&0.0000&8.7530&1.0000\\
$w^{(t)}_4$&0.0000&3.0000&1.2235\\
$w^{(t)}_5$&0.0000&2.0000&4.8117\\
$w^{(t)}_6$&3.9030&1.5466&1.0000\\
$w^{(t)}_7$&0.0000&0.0000&0.0000\\
$w^{(t)}_8$&0.0000&2.0000&4.0000\\
$w^{(t)}_9$&8.0000&0.0000&0.0000\\
$w^{(t)}_{10}$&0.0000&1.0000&0.4531\\
$w^{(t)}_{11}$&0.0000&3.0000&3.0000\\
$w^{(t)}_{12}$&0.0000&1.0000&1.0000\\
$w^{(t)}_{13}$&0.0581&0.0000&0.0000\\
$w^{(t)}_{14}$&3.0000&0.0000&0.0000\\
$w^{(t)}_{15}$&1.0000&1.0000&1.0000\\
$w^{(t)}_{16}$&0.0000&3.0000&6.0000\\
$w^{(t)}_{17}$&0.0000&0.0000&4.0000\\
$w^{(t)}_{18}$&0.6396&5.0000&5.0000\\
\hline
$||\mathbf{w}^{(t)}-\mathbf{w}^{(t-1)}||$&0.0000&0.0066&0.0011\\
\hline
\end{tabular}
\end{table}

We also examined the proposed algorithm for 3 instances Ruspini 75, Bongartz 287 and TSPLIB 654, which can be found in Beasley \cite{B90}. All the weights, costs and upper bounds are randomly generated in the interval $[1,10]$ . The range of coordinates of the points in these instances are given in the third column of Table \ref{medians}. Since the median is in the convex hull of the existing points, therefore, for each instances, the given point $\mathbf{\bar{x}}$ should be chosen in the range of existing points, otherwise the inverse problem is infeasible. Table \ref{medians} also shows the median points and their objective functions which obtained by the Weiszfeld algorithm with tolerance $\epsilon=0.001$. 

\begin{table}\caption{ The coordinate ranges of existing points and median points of instances from \cite{B90} with varying $L_p$ norms and random weights.}\label{medians}
\centering
\small %\scalebox{1}{
\begin{tabular}{cccccc}
\hline
Instance&$n$&Coordinate ranges&p&$\mathbf{x^{(0)}}=(x^{(0)},y^{(0)})$&$f(\mathbf{x}^{(0)},\mathbf{w}^{(0)},\mathbf{A})$\\
\hline
Ruspini&75&$(5,5)\ to\ (110,110)$&2&$(56.7635,102.9328)$&23630.4385\\
&&&3&$(62.3770,103.0802)$&22543.1684\\
&&&5&$(67.2228,103.7886)$&21855.8203\\
&&&8&$(68.9591,104.1744)$&21566.5435\\
\hline
Bongartz&287&$(5,5)$ to $(48,48)$&2&$(22.3657,32.2146)$&13946.0271\\
&&&3&$(22.2903,32.2109)$&13226.8711\\
&&&5&$(22.2175,32.1861)$&12843.64547\\
&&&8&$(22.1757,32.1678)$&12709.5519\\
\hline
TSPLIB&654&$(1000,1000)$ to $(5000,5000)$&2&$(3410.0620,3689.1907)$&9085118.3175\\
&&&3&$(3455.5845,3704.2810)$&8463911.9145\\
&&&5&$(3389.6068,3116.5265)$&8327910.9097\\
&&&8&$(3379.4386,3085.4853)$&8203351.1022\\
\hline
\end{tabular}
\end{table}

Tables \ref{resultl2} to \ref{resultl8} contains the results obtained by Algorithm A1 for varying $L_p$ norms and given points $\mathbf{\bar{x}}$. In these tables the columns with heading $\Delta w^{(t)}$ indicate the norm of difference between the weights in the last two iterations of algorithm, i.e. $\Delta w^{(t)}=||\mathbf{w}^{(t)}-\mathbf{w}^{(t-1)}||$.  Table \ref{resultl2} shows the results those obtained by our algorithm and method of Burkard et al. \cite{BGG10} for $L_2$ norm. Comparing the results indicates that the method of Burkard et al. \cite{BGG10} is faster than our algorithm. However, their method can just be applied for $L_2$ norm.

\begin{table}\caption{ The results for the instances from \cite{B90} for $L_2$ norm.}\label{resultl2}
\centering
\tiny
 %\scalebox{1}{
\begin{tabular}{ccc|ccc|ccccc}
\hline
&&&\multicolumn{3}{c|}{Method of Burkard et al. \cite{BGG10}}&\multicolumn{5}{c}{Algorithm A1}\\ \hline
Instance&$n$&$\mathbf{\bar{x}}=(\bar{x},\bar{y})$&$\mathbf{x_B}$&$C_B$&CPU&$\mathbf{{x}}^{(t)}$&$t$& $C^*$&$\Delta w^{(t)}$&CPU\\
&&&&&(in sec)&&&&&(in sec)\\
\hline
Ruspini&75&$(50,50)$&(49.964, 50.038)&918.572&0.0263&$(50.003, 50.006)$&21&918.125&0.0058&1.0773\\
 & &$(80,20)$&$(80.052, 20.034)$&2124.434&0.0300&(80.054, 20.068)&20&2121.925&0.0048&0.8295\\
 & &$(20,80)$&$(20.011, 80.022)$&1226.341&0.0411&(20.025, 80.014)&23&1225.508&0.0087&1.0798\\
\hline
Bongartz&287&$(15,35)$&$(14.999,34.996)$&5111.552&0.0554&(15.000,34.995)&19&5110.484&0.0000&1.7481\\
& &$(10,40)$&$(9.996,39.993)$&7836.151&0.0605&(10.007,40.004)& 27&7835.434&0.0000&4.1123\\
& &$(30,20)$&$(29.965,20.003)$&7813.039&0.0820&(29.956,20.016)& 20&7811.124&0.0030&2.3278\\
\hline
TSPLIB&654&$(2000,4000)$&$(2000.036,3999.988)$&1668.822&0.1054&(2000.041,3999.944)&25&1668.666&0.0098&10.2197\\
&&$(1500,1500)$&$(1500.192,1500.108)$&19772.922&0.1288&(1500.284,1500.631)&27&19772.592&0.0064&10.3044\\
&&$(3500,3500)$&$(3499.858,3499.994)$&680.945&0.1460&(3500.050,3500.037)&27&680.716&0.0042&10.0377\\
\hline
\end{tabular}
\end{table}

\begin{table}\caption{ The results for the instances from \cite{B90} for $L_3$ norm.}\label{resultl3}
\small
\centering
 %\scalebox{1}{
\begin{tabular}{ccc|ccccc}
\hline
Instance&$n$&$\mathbf{\bar{x}}=(\bar{x},\bar{y})$&$\mathbf{{x}}^{(t)}$&$t$& $C^*$&$\Delta w^{(t)}$&CPU\\
&&&&&&&(in sec)\\
\hline
Ruspini&75&$(50,50)$&(49.9796,50.0447)&21&797.8499&0.0036&0.7683\\
 & &$(80,20)$&(79.9370,19.9753)&20&2095.7514&0.0084&0.8642\\
 & &$(20,80)$&(20.0039,80.0050)&22&989.6491&0.0007&0.7568\\
\hline
Bongartz&287&$(15,35)$&(15.0005,35.0006)&30&5021.4630&0.0000&4.5220\\
& &$(10,40)$&(10.0045, 40.0008)&25&	7812.9340&0.0081&3.8931\\
& &$(25,25)$&(24.9987, 25.0107)&19&	4545.8444&0.005&2.9389\\
\hline
TSPLIB&654&$(2000,4000)$&(2000.0102, 3999.9955)&27&	2573.2310&0.0054&10.7947\\
&&$(1500,1500)$&(1500.1864, 1500.1759)&34&	19698.5450&0.0084&15.5954\\
&&$(3500,3500)$&(3500.0449, 3500.0118)&35&	1236.0234&0.0054&15.7235\\
\hline
\end{tabular}
\end{table}

\begin{table}\caption{ The results for the instances from \cite{B90} for $L_5$ norm.}\label{resultl5}
\centering
\small
 %\scalebox{1}{
\begin{tabular}{ccc|ccccc}
\hline
Instance&$n$&$\mathbf{\bar{x}}=(\bar{x},\bar{y})$&$\mathbf{{x}}^{(t)}$&$t$& $C^*$&$\Delta w^{(t)}$&CPU\\
&&&&&&&(in sec)\\
\hline
Ruspini&75&$(50,50)$&(49.9887, 50.0299)&15&803.0411&0.0064&	0.6311\\
 & &$(80,20)$&(79.8436, 19.9586)&18&	2047.4155&0.0094&1.3337\\
 & &$(20,80)$&(20.0017, 80.0112)&17&	776.3165&0.0078&0.6362\\
\hline
Bongartz&287&$(15,35)$&(15.0004, 35.0002)&30&4858.1739&0.0023&4.3586\\
& &$(10,40)$&(9.9993, 39.9963)&24&	7750.9010&0.0065& 3.6254\\
& &$(25,25)$&(24.9986, 25.0089&18&	4305.7498&0.0039&	1.6300\\
\hline
TSPLIB&654&$(2000,4000)$&(2000.0134, 3999.9885)&21&3791.6483&0.0091&8.2931\\
&&$(1500,3500)$&(1500.0089, 3500.0034)&36&8707.6269&0.0049&20.4489\\
&&(1500, 4000)&(1500.0035, 4000.0007)&35&7868.2359&0.0078&17.2575\\
\hline
\end{tabular}
\end{table}

\begin{table}\caption{ The results for the instances from \cite{B90} for $L_8$ norm.}\label{resultl8}
\centering
\small
 %\scalebox{1}{
\begin{tabular}{ccc|ccccc}
\hline
Instance&$n$&$\mathbf{\bar{x}}=(\bar{x},\bar{y})$&$\mathbf{{x}}^{(t)}$&$t$& $C^*$&$\Delta w^{(t)}$&CPU\\
&&&&&&&(in sec)\\
\hline
Ruspini&75&$(50,50)$&(49.9987, 50.0048)&15&860.9878&0.0022&0.5870\\
 & &$(80,20)$&(79.9531, 20.0733)&18&2025.5144&0.0054&1.5873\\
 & &$(20,80)$&(19.7536, 80.2957)&10&624.4426&0.0098&0.5051\\
\hline
Bongartz&287&$(15,35)$&(15.0013,35.0000)&26&4729.2314&0.0086&3.2491\\
& &$(10,40)$&(9.9986, 39.9981)&26&7687.0110&0.0079&	3.7591\\
& &$(25,25)$&(25.0002, 25.0075)&22&4161.0671&0.0011&2.8497\\
\hline
TSPLIB&654&$(1500,3500)$&(3500.0117, 3499.9914)&40&9191.1963&0.0099&25.3061\\
%&&$(2000,4000)$&\\
%&&$(3500,3500)$&\\
\hline
\end{tabular}
\end{table}

\begin{figure}[h]
\begin{center}
%\vspace{-1cm}
  \begin{tabular}{cc}
  \hspace{-0.5cm}
 \includegraphics[width=7cm]{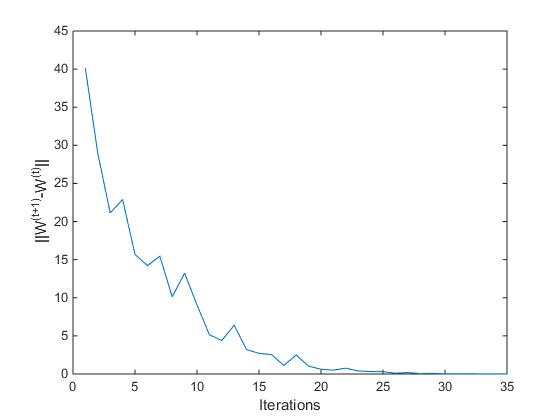}&
 \hspace{-0.5cm}
  \includegraphics[width=7cm]{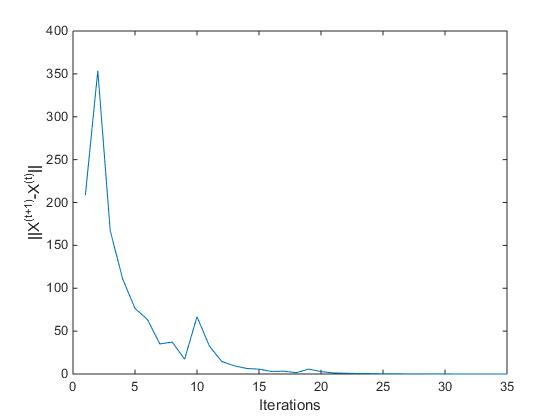}\vspace{-0.2cm}\\ 
  \hspace{-0.5cm}
  (a)& \hspace{-1cm}(b)
  \end{tabular}\caption{The changes of difference weights and obtaining points in two consecutive iterations for instance TSPLIB with $L_3$ norm and  $\mathbf{\bar{x}}=(3500,3500)$.
  (a) The change of difference weights. 
  (b) The change of difference obtaining points.}\label{error}
\end{center}
\end{figure}

As we mentioned in Theorem \ref{stopx}, the stopping condition can be replaced by $||\mathbf{{x}}^{(k)}-\mathbf{{x}}^{(k-1)}||\leq \epsilon$. The results show $\mathbf{{x}}^{(t)}$ is very close to $\mathbf{\bar{x}}$ and in the most cases this alternative stopping condition also holds. To make more visible the comparing of the two stopping conditions, the variations of difference of weights and obtaining points in two consecutive  iterations for instance TSPLIB with $L_3$ norm and  $\mathbf{\bar{x}}=(1500,1500)$ are shown in Figure \ref{error}.   

\section{Summary and conclusion}\label{sec5}

In this paper we presented a row generation algorithm for general case of inverse continuous location problems. The convergence and optimality condition have been proved for the given algorithm. In special cases we considered the inverse minisum and minimax location models and some test problems have been solved by the presented algorithm. The results show the algorithm could efficiently find the optimal solutions. %, and in the rest of cases a solution with reasonable cost are found by the algorithm.   

Note that, the presented algorithm can be applied for solving inverse of other cases of the location problems such as inverse covering problem, inverse line location problem and inverse circle location problem. The algorithm can also be used for solving other kind of inverse location models such as hamming and rectilinear norms.  However, the efficiency of the algorithm depends on the solvability of two problems: 1- finding the optimal solution of problem $(P)$, 2- problem $(P_k)$, in iteration $k$. 

The inverse continuous location problems with variable coordinates can also be solved by the same structure algorithm. However, it will be more harder, since problem $(P_k)$ is not linear even for the inverse minisum location problem ( Baroughi-Bonab et al. \cite{BBA10} showed this problem with rectilinear and Chebyshev norms is NP-hard). This case of inverse location models can be considered in the future works.


\begin{thebibliography}{10}

\bibitem{AB11} Alizadeh B., Burkard R.E., Combinatorial algorithms for inverse absolute and vertex 1-center location problems on trees, Networks,  58 (2011), 190-200.

\bibitem{BBA10}Baroughi-Bonab F., Burkard R. E., Alizadeh B., Inverse median location problems with variable coordinates, Central European Journal of Operations Research, 18 (2010), 365-381.

\bibitem{B90} Beasley J.E., OR-Library: distributing test problems by electronic mail. Journal of the Operational Research Society, 41 (1990), 1069-1072.

\bibitem{BGG10} Burkard R. E., Galavii M., Gassner E., The inverse Fermat–Weber problem, European Journal of Operation Research, 206 (2010), 11-17.

\bibitem{BPZ04} Burkard R. E., Pleschiutschnig C., Zhang J. Z., Inverse median problems, Discrete Optimization, 1 (2004), 23-39.

\bibitem{CYZ99} Cai M.C., Yang X.G., Zhang J., The complexity analysis of the inverse center location problem, Journal of Global Optimization, 15 (1999), 213-218.

\bibitem{EM11} Eiselt H. A., Marianov V., Foundations of Location Analysis, Springer, New York, (2011).

\bibitem{F09} Fathali J., Zaferanieh M. and Nezakati A., A
BSSS algorithm for the location problem with minimum square error, {\em Advances in Operations Research}, Volume 2009 (2011), Article ID 212040, 10 pages.  

\bibitem{G10} Galavii M., The inverse 1-median problem on a tree and on a path, Electronic Notes in Discrete Mathematics, 36 (2010), 1241-1248.

\bibitem{GZ12} Guan X. C., Zhang B. W., Inverse 1-median problem on trees under weighted Hamming distance, Journal of Global Optimization, 54 (2012), 75-82.

\bibitem{LMV88} Love R.F., Morris J.G., Wesolowsky G.O., Facilities Location: Models and Methods, North-Holland, New York, (1988).

\bibitem{NFNV18}Nazari M., Fathali J., Nazari M., Varedi-Koulaei S.M., Inverse of backup 2-median problems with variable edge lengths and vertex weight on trees and variable coordinates on the plane, Production and Operations Management, 9 (2018), 115-137.

\bibitem{NF18}Nazari M., Fathali J., Reverse backup 2-median problem with variable coordinate of vertices, Journal of Operational Research and Its Applications 15 (2018), 63-88.

\bibitem{NNNLP19} Nguyen K. T., Nguyen T. H., Nguyen-Thu H., Le T. T., Pham V. H., On some inverse 1-center location problems, 68 (2019),  999-1015.


%\bibitem{ABP09} Alizadeh B., Burkard R.E., Pferschy U., Inverse 1-center location problems with edge length augmentation on trees, Computing,  86 (2009), 331-343.

 
%\bibitem{BPZ08} Burkard R. E., Pleschiutschnig C., Zhang J. Z., The inverse 1-median problem on a cycle, Discrete Optimization, 5 (2007), 242-253.

%\bibitem{N16} Nguyen K.T., Reverse 1-center problem on weighted trees, Optimization, 65 (2016), 253-264.

\bibitem{N16-2} Nguyen K.T., Inverse 1-median problem on block graphs with variable vertex weights, Journal of Optimization Theory and Applications, 168 (2016), 944-957.

\bibitem{NS15} Nguyen K.T., Sepasian A.R., The inverse 1-center problem on trees with variable edge lengths under Chebyshev norm and Hamming distance, Journal of Combinatorial Optimization
32, (2016), 872-884.

\bibitem{OF20} Omidi S., Fathali J., Nazari M., Inverse and reverse balanced facility location problems with variable edge lengths on trees, OPSEARCH, 57, (2020), 261-273.

\bibitem{SR15} Sepasian, A.R., Rahbarnia, F., An O(nlogn) algorithm for the inverse 1-median problem on trees with variable vertex weights and edge reductions. Optimization, 64 (2015), 595-602. 

\end{thebibliography}
\end{document}